# NONPARAMETRIC CHECKS FOR SINGLE-INDEX MODELS

By Winfried Stute and Li-Xing Zhu[1]

*University of Giessen and University of Hong Kong and
East China Normal University at Shanghai*

In this paper we study goodness-of-fit testing of single-index models. The large sample behavior of certain score-type test statistics is investigated. As a by-product, we obtain asymptotically distribution-free maximin tests for a large class of local alternatives. Furthermore, characteristic function based goodness-of-fit tests are proposed which are omnibus and able to detect peak alternatives. Simulation results indicate that the approximation through the limit distribution is acceptable already for moderate sample sizes. Applications to two real data sets are illustrated.

**1. Introduction.** Suppose that a response variable $Y$ depends on a vector $X = (x_1, \ldots, x_p)^T$ of covariates, where $T$ denotes transposition. We may then decompose $Y$ into a function $m(X)$ of $X$ and a noise variable $\varepsilon$, which is orthogonal to $X$, that is, for the conditional expectation of $\varepsilon$ given $X$ we have $\mathbb{E}(\varepsilon|X) = 0$. When $Y$ is unknown, the optimal predictor of $Y$ given $X = \mathbf{x}$ equals $m(\mathbf{x})$. Since in practice the regression function $m$ is unknown, statistical inference about $m$ is an important issue. In a purely parametric framework, $m$ is completely specified up to a parameter. For example, in linear regression $m(\mathbf{x}) = \beta^T \mathbf{x}$, where $\beta$ is an unknown $p$-vector which needs to be estimated from the available data. Slightly more generally we may consider $m(\mathbf{x}) = \Phi(\beta^T \mathbf{x})$, where the link-function $\Phi$ may be nonlinear but is again specified. This is the so-called generalized linear model.

When $\Phi$ remains unspecified, we arrive at a semiparametric model which is more flexible on the one hand and, on the other hand, avoids the curse of dimensionality one faces in fully nonparametric models. The estimator of $\beta$, as well as of the link function $\Phi$, in this so-called single-index model

Received February 2003; revised April 2004.
[1]Supported by a grant of the Chinese Academy of Sciences and Grant HKU7129/00P of the Research Grants Council of Hong Kong.

was studied by among others, Li and Duan [25], Härdle, Hall and Ichimura [16], Ichimura [23] and Hristache, Juditsky and Spokoiny [22]. Related work is [6] and [20]. Clearly, any statistical analysis within the model, to avoid wrong conclusions, should be accompanied by a check of whether the model is valid at all. For the single-index model the diagnostic methods are less elaborate. We only mention Fan and Li [14], Aït-Sahalia, Bickel and Stoker [1] and Xia, Li, Tong and Zhang [38] here but come back to them later. See Discussion 2.6, when we are prepared to compare their approaches and results with ours. The paper by Härdle, Mammen and Proença [19] considers a parametric link structure and therefore does not fall into the area studied in this paper.

In the present paper, we aim at developing some formal tests for model checking when the link function remains unspecified.

For more specified regression models the literature is much more elaborate. To review only a few contributions, Cox, Koh, Wahba and Yandell [8] introduced tests of the null hypothesis that a regression function has a particular parametric structure. Azzalini, Bowman and Härdle [3] considered nonparametric regression as an aid to model checking. Cox and Koh [7] developed spline-based tests of model adequacy. Eubank and Spiegelman [11] considered spline approaches to testing the goodness of fit of a linear model. Simonoff and Tsai [28] proposed diagnostic methods for assessing the influence of individual data values on goodness-of-fit tests based on nonparametric regression. Gu [15] used spline methods in a diagnostic approach to model fitting. Azzalini and Bowman [2] used nonparametric regression to check linear relationships. Eubank and LaRiccia [10] derived properties of two-sided tests in nonparametric regression based on Fourier methods. Härdle and Mammen [17] considered comparisons between parametric and nonparametric fits and used the wild bootstrap for the computation of critical regions. Härdle, Mammen and Müller [18] investigated testing for parametric versus semiparametric modeling in generalized linear models, again using the wild bootstrap.

Note, however, that any test using a nonparametric regression estimator runs into an ill-posed problem requiring the choice of a smoothing parameter. Therefore, an alternative approach was developed which circumvents these problems. To name only a few papers, Bierens [4] proposed to check a parametric regression model by investigating the sum of properly weighted residuals. See also [5] for an informative discussion of the resulting tests when local alternatives are considered. In Stute [33] a method was studied which is based on the integrated regression function and which corresponds to cumulative quantities such as empirical distribution functions or ranks known from other areas in statistics. In this setup the author was able to derive a principle components decomposition of the underlying test process, which is extremely useful for design of optimal tests versus local alternatives and for understanding the impact of the design distribution and the noise variance on the power of the tests. In particular, optimal Neyman–Pearson tests which are based on linear rather than quadratic test statistics can be obtained from this decomposition. Stute, González Manteiga and Presedo Quindimil [35] studied the quality of the distributional approximation of an associated cusum process via the wild bootstrap, while Stute, Thies and Zhu [36] proposed an innovation process approach so as to obtain asymptotically distribution-free and optimal tests. Finally, Stute and Zhu [37] developed nonparametric testing for the validity of a generalized linear model, which is based on a proper transformation of a residual empirical process and which perfectly adapts to a situation when the design vector is elliptically contoured.

In the framework of the single-index model the link function is unknown and, as part of the testing procedure, needs to be estimated in a nonparametric way. From our preceding remarks on ill-posedness, one might conclude that nonparametric estimation of the link function necessarily excludes the possibility of constructing tests which have optimal power versus local alternatives converging to the null model at the rate $n^{-1/2}$. Fortunately, as this paper will show, this pessimistic view is not justified. To obtain such tests, rather than comparing the estimator of $\Phi$ with the hypothetical semiparametric model, we embed the residuals into a cusum process. This summation has a smoothing effect so that our test is much less sensitive than usual to a wrong choice of the bandwidth. At the same time, each residual is properly weighted by a function of the design vector. Our main result, Theorem 2.1, is formulated for a given fixed weight function. Such an approach has a long tradition in statistics. Typically, score tests are first analyzed (and optimized) when the direction from which the alternative tends to the null model is specified. Classical examples are linear one- and two-sample rank statistics or rank correlation statistics. Also, robust tests focussing on a neighborhood of a given family of distributions are designed in this spirit.

Theorem 2.1 not only provides the asymptotic normality of a large class of score statistics, but also yields (up to a remainder) a representation as a sum of i.i.d. variables. From this, when the alternative is specified, we shall be able to choose the weights so as to optimize local power. This discussion will give us a clue as to how to proceed if the alternative model has arbitrary but finite codimension $d$. In such a situation we propose and study a test which is asymptotically distribution-free and shown to be maximin (Corollary 2.2). Since $d$ is arbitrary, Corollary 2.2 covers most situations arising in practice. The i.i.d. representation is also useful for implementation of a proper bootstrap approximation. See Section 3 for some details.

For those readers who prefer omnibus tests, we also discuss (Theorem 2.3) a situation where the deviation from the null model is completely nonparametric. Also, in this case, the local asymptotic power can be derived.

Finally, we include a discussion of how our test behaves when local peak alternatives are to be detected.

The paper is organized as follows. In Section 2 we introduce the basic test statistics and formulate our main results. In Section 3 we report on some simulation results and apply our method to two data sets. Proofs of theoretical results are postponed to Section 4. Readers who want to skip the technical part may consult Section 2 for an informal discussion and some background information on proofs.

**2. Main theorems.** Throughout the paper we assume that the available data $(X_i, Y_i)$, $1 \leq i \leq n$, are independent and have the same distribution as $(X, Y)$. Under the null hypothesis, that is, under the single-index model,

$$Y = \Phi(\beta^T X) + \varepsilon, \tag{2.1}$$

where $\beta$ is an unknown $p$-vector and $\Phi$ is an unspecified link function defined on the real line. The noise variable $\varepsilon$ satisfies

$$\mathbb{E}(\varepsilon|X) = \mathbb{E}(\varepsilon|\beta^T X) = 0, \tag{2.2}$$

which is tantamount to saying that

$$\mathbb{E}(Y|X) = \mathbb{E}(Y|\beta^T X) = \Phi(\beta^T X). \tag{2.3}$$

Note that (2.2) allows $\varepsilon$ to depend on $X$ so that (2.1) may include heteroscedastic errors. The first equation in (2.3) features the projection pursuit character of the single-index model in that the conditional mean of $Y$ given $X$ only depends on a proper projection of $X$.

To motivate our approach, assume for a moment that we already have an estimator $\hat{\beta}$ of $\beta$. Replacing $\beta^T X_i$ with $\hat{\beta}^T X_i$, we could try to estimate $\Phi$ through a Nadaraya–Watson estimator $\hat{\Phi}$ or a local linear smoother as discussed, for example, in [13]. The disadvantage of these smoothers, at least in our context, comes from the fact that the distribution of $\hat{\beta}$, as well as $X$, will likely have an effect on the distribution of our test statistic, even in the limit. This phenomenon is well known in many other statistical problems, when unknown parameters need to be estimated. Typically, the effect on the distributional character requires some correction through a proper transformation of the test statistic. See, for example, [34]. Moreover, the ratio structure of these estimators $\hat{\Phi}$ creates some technical problems when the denominator is small, that is, when **x** lies in a region of low density. From time to time some structural assumptions on level sets are imposed, but when it comes down to estimation, these assumptions can hardly be justified for $\hat{\Phi}$. To avoid all these nasty side effects, we decided to use an estimator of $\Phi$ which employs a transformation of $\hat{\beta}^T X_i$ to a variable which is approximately uniform on the unit interval $(0,1)$. In other words, we incorporate a transformation which makes everything distribution-free, as far as the distribution of $\beta^T X$ is concerned. This estimator is a symmetrized nearest-neighbor (NN) estimator. Its consistency was proved by Yang [39], while Stute [32] provided the asymptotic normality. In these papers, the regression function itself was, of course, the target and the distribution-freeness only applies to the random deviation but not to the bias term. In the context of the present paper, $\hat{\Phi}$ only appears as a tool to define the residuals. When we consider a properly weighted sum of the residuals, averaging yields a smaller variance to the effect that we may choose smoothing parameters so that at the same time the bias becomes negligible and the variance part remains as the only nonnegligible source of error. This more or less enables us to construct tests which have nontrivial power when the alternatives approach the null model at the rate $n^{-1/2}$.

To motivate our approach on a more technical level, assume that $\beta^T X$ has a continuous distribution function $F^\beta$, that is,

$$F(x) \equiv F^\beta(x) := \mathbb{P}(\beta^T X \leq x), \qquad x \in \mathbb{R}.$$

Here $\mathbb{P}$ denotes a probability measure defined on a space $(\Omega, \mathcal{A})$ carrying all random variables which may appear. Denote by $F^{-1}$ the quantile function of $F$:

$$F^{-1}(u) = \inf\{x \in \mathbb{R} : F(x) \geq u\}, \qquad 0 < u < 1.$$

Put $U := F(\beta^T X)$. By continuity of $F$, the variable $U$ has a uniform distribution on $(0, 1)$. Setting

$$\psi = \Phi \circ F^{-1},$$

equation (2.1) becomes (with probability one)

$$Y = \psi(U) + \varepsilon.$$

In terms of regression, this may be expressed as

$$m(\mathbf{x}) \equiv \mathbb{E}(Y|X = \mathbf{x}) = \Phi(\beta^T \mathbf{x}) = \psi(u),$$

where

$$u = F(\beta^T \mathbf{x}) \text{ and } \psi(u) = \mathbb{E}(Y|F(\beta^T X) = u).$$

Therefore, the kernel estimator for $\psi$ at $0 < u < 1$ becomes

$$\hat{\psi}_n(u) = \frac{1}{n} \sum_{i=1}^n Y_i K_h(u - U_i),$$

where

$$K_h(v) = \frac{1}{h} K\left(\frac{v}{h}\right)$$

and $K$ is a symmetric kernel on the real line integrating to one, while $h = h_n > 0$ is a bandwidth. The random variables

$$U_i = F^\beta(\beta^T X_i)$$

are i.i.d. from the uniform distribution on $(0,1)$. Since $F^\beta$ and $\beta$ are unknown, $\hat\psi_n$ cannot be our final estimator. For this, replace $\beta$ by some estimator $\hat\beta$ and $F = F^\beta$ by the empirical distribution function $F_n$ of $\hat\beta^T X_i$, $1 \leq i \leq n$. This yields

$$\hat U_i := F_n(\hat\beta^T X_i), \qquad 1 \leq i \leq n,$$

with corresponding estimator

$$\psi_n(u) = \frac{1}{n} \sum_{i=1}^n Y_i K_h(u - \hat U_i).$$

This estimator is related to that in [32], up to the fact that there univariate $X_i$'s were considered and no preliminary projection was required. The $\hat U_i$'s are the normalized ranks pertaining to the projected values $\hat\beta^T X_i$. Since these values depend on the random $\hat\beta$, existing results on rank statistics cannot give us easy access to the analysis of our final test statistic, in particular, since the $\hat U_i$'s appear as part of the smoothed function $\psi_n$ at $u$.

Worse than that, we have to evaluate $\psi_n$ at each $\hat U_j$. This finally leads to the residuals

$$\hat\varepsilon_j = Y_j - \psi_n(\hat U_j), \qquad 1 \leq j \leq n.$$

Actually, to reduce a possible bias, we shall consider estimators $\psi_n^{(j)}$ computed in the same way as $\psi_n$, but with the $j$th datum deleted from the observations. Hence, the residuals are to be redefined as

$$\hat\varepsilon_j = Y_j - \psi_n^{(j)}(\hat U_j), \qquad 1 \leq j \leq n.$$

The mathematical analysis of $\psi_n^{(j)}(\hat U_j)$ and, hence, of $\hat\varepsilon_j$ requires careful study of the local properties of $F_n$ evaluated at $\hat\beta^T X_i$. The oscillation behavior for the ordinary empirical process has been investigated in detail in [30, 32]. In the present situation we need to study the fluctuations of empirical measures over halfspaces rather than quadrants.

Our final test statistic will be of the form

$$\hat T_n = n^{-1/2} \sum_{j=1}^n \hat\varepsilon_j W_j.$$

The weights $W_j$ will be of the form $W_j = W(X_j)$. The function $W$ is a smooth function defined on $\mathbb{R}^p$. A discussion of how to choose $W$ in a testing situation is postponed to the end of this section. Under the null model (2.2), we may expect that $\hat{T}_n$ behaves similarly to

$$T_n = n^{-1/2} \sum_{j=1}^{n} \varepsilon_j W_j.$$

Since $W_j$ is orthogonal to $\varepsilon_j$, $T_n$ is centered. Hence, we may expect that also $\hat{T}_n$ fluctuates around zero under (2.2). Under (local) alternatives, the $\hat{\varepsilon}_j$ also comprise quantities which hopefully are not orthogonal to the $W_j$'s. If we choose $W$ in a proper way, this fact will guarantee nontrivial power of the test.

More specifically, we shall first consider models of the type

(2.4) $$Y_{in} = \Phi(\beta^T X_i) + n^{-1/2} s(X_i) + \varepsilon_i, \qquad 1 \leq i \leq n,$$

where the $(X_i, \varepsilon_i)$ are i.i.d. satisfying

(2.5) $$\mathbb{E}(\varepsilon_i | X_i) = 0 \qquad \text{for } 1 \leq i \leq n.$$

The function $\Phi$, as well as the parameter $\beta$, remain unspecified, as will be the distribution of $X_i$ and $\varepsilon_i$. The function $s$ may or may not be specified. When $s \equiv 0$, the single-index model holds. For specified alternatives, we shall later discuss how to choose $W$ in order to maximize local power.

So far we have not discussed how to estimate $\beta$. We shall come back to this point in Section 3 when we apply our method in a simulation study and to real data. In fact, the discussion of $\hat{\beta}$ may be delayed since our assumptions on $\hat{\beta}$ are very general and do not assume any particular form for $\hat{\beta}$.

We now state the assumptions needed for Theorem 2.1 below. For this, put, for $0 < u < 1$,

$$\bar{W}(u) = \mathbb{E}[W(X)|U = u], \qquad \bar{s}(u) = \mathbb{E}[s(X)|U = u].$$

THEOREM 2.1. *Assume that* (2.4), (2.5) *and the following conditions hold:*

A     (i) *$\psi, \bar{s}$ and $\bar{W}$ are twice continuously differentiable.*
       (ii) *$YW(X)$ and $\varepsilon W(X)$ have finite second moments.*
B     (i) *$\mathbb{E}\|X\|^\gamma < \infty$ for some $\gamma > 2$.*
       (ii) *For all $\theta$ in a neighborhood of $\beta$, the variables $\theta^T X$ have continuous densities $f^\theta$ which are uniformly bounded.*
       (iii) *The distribution functions $F^\theta$ of $\theta^T X$ are continuous in $\theta$ at $\theta = \beta$.*
       (iv) *The estimator $\hat{\beta}$ satisfies $n^{1/2}(\hat{\beta} - \beta) = O_\mathbb{P}(1)$.*
C     (i) *$n^{1/2} h^2 \to 0$ and $h^{-1} n^{-1/2 + 1/\gamma} \to 0$.*
       (ii) *$K$ is a symmetric kernel with compact support, twice continuously differentiable with $\int K = 1$. Furthermore, $K$ is nonincreasing on the positive real numbers.*

*Then we have*

(2.6) $$\hat{T}_n = \mu + n^{-1/2} \sum_{i=1}^n \varepsilon_i [W_i - \bar{W}(U_i)] + o_\mathbb{P}(1)$$

*and, therefore, by the CLT,*

$$\hat{T}_n \to \mathcal{N}(\mu, \sigma^2) \quad \text{in distribution,}$$

*where*

$$\sigma^2 = \mathbb{E}\{\varepsilon^2 [W(X) - \bar{W}(U)]^2\}$$

*and*

$$\mu = \mathbb{E}\{[s(X) - \mathbb{E}(s(X)|U)]W(X)\}.$$

A discussion of A–C will be postponed until the end of this section.

The drift comprises the deviation of $s(X)$ from the space of variables spanned by $\beta^T X$. Under the single-index model, the bracket equals zero and so does $\mu$. Also, $W(X)$ should not depend on $X$ through $\beta^T X$, since then also $\mu = 0$. The variance does not depend on $s$ but, among other things, measures the deviations between $W(X_j)$ and the projected values $\bar{W}(U_j)$. The limit variance $\sigma^2$ also does not depend on the unknown $\Phi$. A consistent estimator of $\sigma^2$ is obtained by

$$\sigma_n^2 = \frac{1}{n} \sum_{j=1}^n \hat{\varepsilon}_j^2 [W(X_j) - \bar{W}_n^{(j)}(\hat{U}_j)]^2,$$

where $\bar{W}_n^{(j)}$ is defined similarly to $\psi_n^{(j)}$. Just replace $Y_i$ with $W(X_i)$ in the definition of the NN-estimator. Putting

$$\bar{T}_n := \hat{T}_n / \sigma_n,$$

we then obtain

$$\bar{T}_n \to \mathcal{N}(C, 1) \quad \text{in distribution,}$$

with

$$C = \mu/\sigma.$$

The null model is rejected at level $\alpha$ if

$$|\bar{T}_n| \geq \lambda_{1-\alpha/2} \equiv \lambda,$$

where $\lambda$ is the $(1 - \frac{\alpha}{2})$-quantile of the standard normal distribution function $\Phi$. Hence, the asymptotic power of $|\bar{T}_n|$ against the local alternatives (2.4) equals $1 - [\Phi(C + \lambda) - \Phi(C - \lambda)]$. This is a monotone function of $|C|$. Thus, we should select the weight function $W$ in a way that makes $C^2$ as large as possible. If we write, in an obvious notation,

$$C^2 = C^2(s, W) = \frac{\mu^2(s, W)}{\sigma^2(W)},$$

it is easy to determine the optimal solution of our problem when the $\varepsilon$'s are independent of $X$, that is, if the homoscedastic case holds. Then the above ratio equals

$$\frac{\mu^2(s, W)}{\mathbb{E}\varepsilon^2 \mathbb{E}[W(X) - \bar{W}(U)]^2},$$

and the Cauchy–Schwarz inequality immediately yields that the optimal weight function $W_0$ equals, up to a constant factor, the function $s$:

(2.7) $$W_0(\mathbf{x}) = s(\mathbf{x}).$$

Next we study an important extension of (2.4). For this, let $s_1, \ldots, s_d$ be any finite number of functions, where $d \geq 1$. In applications, these functions may constitute a possible (mean) dependence of $Y$ on $X = \mathbf{x}$ other than projections of $\mathbf{x}$. For example, some of the $s$-functions may be quadratic forms, and others may be in charge of possible interactions between coordinates of $X$.

Instead of (2.4), we therefore consider the more complex model

(2.8) $$Y_{in} = \Phi(\beta^T X_i) + n^{-1/2} \sum_{j=1}^{d} \gamma_j s_j(X_i) + \varepsilon_i, \qquad 1 \leq i \leq n,$$

where $\beta \in \mathbb{R}^p$, $\gamma_1, \ldots, \gamma_d \in \mathbb{R}$ are unknown parameters and $\Phi$ is a nonspecified link function. The null model thus corresponds to

$$H_0 : \gamma_1 = \cdots = \gamma_d = 0.$$

In the following we shall derive maximin tests for $H_0$ versus $\|\gamma\| \geq c$, where $\|\cdot\|$ is a proper norm and $\gamma^T = (\gamma_1, \ldots, \gamma_d)$. Needless to say, such test problems have been well studied in the context of linear regression. The present situation is much more complex since now the null model is the semiparametric single-index model. To the best of our knowledge, the following setup provides the first maximin-test in semiparametric regression. For this, and in view of (2.7), we consider the score-statistics $\hat{T}_n^j$ pertaining to $W = s_j$, $j = 1, \ldots, d$. Put

$$\hat{T}_n = (\hat{T}_n^1, \ldots, \hat{T}_n^d)^T.$$

Theorem 2.1 implies that, under (2.8) (in the homoscedastic case), we have in distribution, as $n \to \infty$,

(2.9) $$\hat{T}_n \to \Sigma \begin{pmatrix} \gamma_1 \\ \vdots \\ \gamma_d \end{pmatrix} + \mathcal{N}_d(\mathbf{0}, \rho^2 \Sigma).$$

Here, $\Sigma = (\sigma_{ij})_{1 \leq i,j \leq d}$ with

$$\sigma_{ij} = \mathbb{E}\{[s_i(X) - \mathbb{E}(s_i(X)|U)][s_j(X) - \mathbb{E}(s_j(X)|U)]\},$$

$\mathcal{N}_d$ denotes a normal distribution on $\mathbb{R}^d$ and $\rho^2 = \mathbb{E}\varepsilon^2$. Assertion (2.9) exhibits that, in the limit, $\hat{T}_n$ is a standard Gaussian shift model. Distributional characteristics of the model (2.8) only appear through the (estimable) covariance matrix. This observation once again supports our approach, in particular, the use of the NN-smoother and the rank transformation.

We may now use existing maximin-theory to obtain optimal tests for $H_0$. See, for example, [29], Theorem 30.2. For this define $\sum_n = (\sigma_{ijn})_{1 \leq i,j \leq d}$ through

$$\sigma_{ijn}^2 = \frac{1}{n} \sum_{k=1}^{n} \hat{\varepsilon}_k^2 [s_i(X_k) - \bar{s}_i^{(k)}(\hat{U}_k)][s_j(X_k) - \bar{s}_j^{(k)}(\hat{U}_k)].$$

COROLLARY 2.2. *For a given significance level $0 < \alpha < 1$, the test*

$$t = \mathbb{1}_{\{\hat{T}_n^T \sum_n^{-1} \hat{T}_n \geq c_\alpha\}}$$

*is a maximin $\alpha$-test for $H_0$ versus $H_1 : \gamma^T \Sigma \gamma \geq \rho^2 a$. Here $c_\alpha$ is the $(1-\alpha)$-quantile of the chi-square random variable $\chi_d^2$ with d degrees of freedom. The asymptotic maximin power is given by $\mathbb{P}(\chi_d^2(a) \geq c_\alpha)$, where now a is the noncentrality parameter.*

Since the codimension $d$ is arbitrary, Corollary 2.2 covers many examples of interest. Some, for example, interaction alternatives, are studied in Section 3. For those who prefer omnibus tests, we now discuss a class of tests which has reasonable power over a nonparametric class of alternatives.

Hence, we come back to (2.4) but leave $s$ unspecified. In order to achieve power, we need to consider a family of weight functions $\{W_\gamma\}_\gamma$ guaranteeing that at least one $W_\gamma$ is able to detect a possible deviation of $s(X) - \bar{s}(U)$ from zero. A class of (smooth) score functions which has found a lot of interest in classical empirical process theory is the family of trigonometric functions. This led to an intensive study of the empirical characteristic function. See, for example, [12] for a nice review and further applications. In our context, $W_\gamma$ therefore becomes

$$(2.10) \qquad W(\gamma, \mathbf{x}) = \exp[i\gamma^T \mathbf{x}],$$

where $i$ is the complex unit and $\gamma \in \mathbb{R}^p$. If we take only finitely many $\gamma$'s, we may conceive, as in Corollary 2.2, asymptotically distribution free $\chi^2$-tests. To handle a nonparametric alternative, we have to let $\gamma$ vary over $\mathbb{R}^p$. Hence, we come up with a stochastic process

$$\hat{T}_n(\gamma) := n^{-1/2} \sum_{j=1}^{n} \hat{\varepsilon}_j W_j(\gamma),$$

where $W_j(\gamma) = W(\gamma, X_j)$. Note that $\hat{T}_n$ has continuous sample paths in $\gamma$. The convergence of the finite-dimensional distributions again follows from (2.6). Tightness is not difficult as long as $\gamma$ varies in a compact set, since the $W(\gamma, \mathbf{x})$ are smooth functions in $\gamma$ and $\mathbf{x}$. For detailed arguments, one needs to check the proof of Theorem 2.1 and show that the remainders are uniformly small on compact $\gamma$-sets, while the leading terms are uniformly continuous. After all this we then come up with the following result.

THEOREM 2.3. *Under the assumptions of Theorem 2.1, the stochastic processes $\{\hat{T}_n(\gamma) : \gamma \in \mathbb{R}^p\}$ converge in distribution (on compact sets) to a continuous Gaussian stochastic process $\hat{T}_\infty$ such that*

$$(2.11) \qquad \mu(\gamma) \equiv \mathbb{E}\hat{T}_\infty(\gamma) = \mathbb{E}\{[s(X) - \bar{s}(U)]W(\gamma, X)\}$$

*and*

$$\mathrm{Cov}(\hat{T}_\infty(\gamma_1), \hat{T}_\infty(\gamma_2)) = \mathbb{E}\{\varepsilon^2[W(\gamma_1, X) - \bar{W}(\gamma_1, U)][W(\gamma_2, X) - \bar{W}(\gamma_2, U)]\}.$$

A Kolmogorov–Smirnov (KS) type test rejects $H_0$ if

$$\tilde{T}_n \equiv \sup_\gamma |\hat{T}_n(\gamma)| \geq c_\alpha,$$

where $c_\alpha$ is the $(1-\alpha)$-quantile of $\sup_\gamma |\hat{T}_\infty(\gamma)|$ under $H_0$, that is, $s \equiv 0$. Since this test is no longer distribution-free, a bootstrap approximation is recommended. See Section 3 for further details. For power considerations, we expand $\mu(\gamma)$ at $\beta$ yielding

$$\mu(\gamma) = \mathbb{E}\{[s(X) - \bar{s}(U)]W(\beta, X)\exp[i(\gamma - \beta)^T X]\}$$
$$\sim \mathbb{E}\{[s(X) - \bar{s}(U)]W(\beta, X)\}$$
$$+ i(\gamma - \beta)^T \mathbb{E}\{(s(X) - \bar{s}(U))W(\beta, X)X\}.$$

The first integral vanishes, since $s(X) - \bar{s}(U)$ is orthogonal to the space of random variables measurable w.r.t. $\beta^T X$. The second (vector-valued) integral $I = I(s)$, say, usually does not vanish so that, for example,

$$\sup_\gamma |\mu(\gamma)| \sim \sup_\gamma \|\gamma - \beta\| \|I\| > 0.$$

This property guarantees that the KS-test has asymptotic power $> \alpha$ uniformly for all $s$ for which $\|I(s)\|$ is bounded away from zero.

Needless to say, a version of Theorem 2.3 also holds for other parametric families of functions $W(\gamma, \cdot)$. We focussed on trigonometric functions since they are at the same time smooth and measure determining and allow for a simple expansion of the drift function.

Though our results cover a large class of local alternatives, people sometimes are interested in detecting so-called "peak alternatives." For this, one needs to consider shift functions $s$ which depend on $n$ in such a way that, as $n \to \infty$, $s_n$ (weakly) converges to a Dirac function or a linear combination of such functions. A typical candidate is

$$s_n^0(\mathbf{x}) = a_n^{-p} \varphi\left(\frac{\mathbf{x} - \mathbf{x}_0}{a_n}\right), \tag{2.12}$$

where $a_n \to 0$ but $na_n^p \to \infty$. The "density" $\varphi$, as well as $\mathbf{x}_0$, the center of the peak, remain unspecified. The test process $\hat{T}_n(\cdot)$ may also serve as a basis to detect alternatives (2.8), where some of the $s_j$'s are of "global type," that is, do not depend on $n$. Others may be of type (2.12). Since the covariance is not affected by the shift, the limit covariance remains the same as in Theorem 2.3. Relevant proofs only deal with the null model so that no changes are required. The shift only enters into Lemmas 4.4 and 4.5, resulting in Corollary 4.6. Taking into account the local flavor of (2.12), these lemmas need some minor modifications resulting, under $s = s_n^0$ from (2.12), in the drift function

$$\mu(\gamma) = \mathbb{E}\hat{T}_\infty(\gamma) = [s(\mathbf{x}_0) - \bar{s}(u_0)]W(\gamma, \mathbf{x}_0)\varphi(\mathbf{0})f(\mathbf{x}_0), \tag{2.13}$$

where $f$ is the density of $X$. Here $u_0 = F(\beta^T \mathbf{x}_0)$. Details are omitted. The function (2.13) nicely features the components which determine the power of the test when $s$ equals (2.12):

- The $X$-density at $\mathbf{x}_0 : f(\mathbf{x}_0)$.
- The "height" of the peak at $\mathbf{x}_0 : \varphi(\mathbf{0})$.
- The deviation of $s$ from the null model at $\mathbf{x}_0 : s(\mathbf{x}_0) - \bar{s}(u_0)$.

If we let $\gamma$ vary over a large compact set, the Kolmogorov–Smirnov test associated with $\hat{T}_n$ is able to detect peak alternatives which converge to the null model at the rate $n^{-1/2}$. The asymptotic power exceeds $\alpha$ but is less than one, depending on the three components discussed above. In particular, our approach yields the correct asymptotics. This finding should be compared with other approaches, where, for much simpler purely parametric regression models, alternatives had to converge to the null model at a rate lower than $n^{-1/2}$. See, for example, [21] and references therein. Not unexpectedly, the power then converges to one.

We continue with some comments on A–C.

REMARK 2.4. Condition A comprises standard smoothness and moment assumptions on the involved functions. Condition B requires some weak conditions on the design vector and on $\hat{\beta}$. In C, $\sqrt{n}h^2 \to 0$ will be needed to make the bias tend to zero. The second assumption on $h$ will be needed to control the fluctuations of the random sums. In view of the fact that we always deal with standardized sums and also that large $X_i$'s may enter the statistics, some connection with the tails of $X$ (in terms of $\gamma$) are natural.

The conditions on $K$ are also standard. The monotonicity of $K$ guarantees that $K'$ has identical signs on the positive and negative reals. Moreover, $K'(0) = 0$. In other words, $K$ may be decomposed into two parts, each of which is compactly supported, by the positive and negative real lines, respectively, and having identical signs there. This property is useful in proofs when, after Taylor's expansion, $K'$ appears as a smoothing kernel.

REMARK 2.5. The conditions on $h$ are weak and are satisfied for a large class of bandwidths. A referee pointed out that this fact could be interpreted as a kind of robustness of the method w.r.t. the choice of $h$. In particular, they do not depend, as in related work, on the dimension $p$ of the $X$-vector or higher degrees of smoothness of the involved functions. We may choose $h$ so that $n^{1/2}h^2$ and $h^{-1}n^{-1/2+1/\gamma}$ are of the same order. This yields

$$h \sim n^{-1/3+1/3\gamma}.$$

In the next section we propose two adaptive methods of bandwidth choice which worked very well in our simulation study. If we are not only interested in maximizing power for a given alternative, we may choose a $W$ with compact support. In this way the test is robust against outliers among the $X_i$'s. Our proof then works with $\gamma = \infty$, that is, $\frac{1}{\gamma} = 0$. In this case, $h \sim n^{-1/3}$.

DISCUSSION 2.6. It is time to compare our approach and results with those of Fan and Li [14], Aït-Sahalia, Bickel and Stoker [1] and Xia, Li, Tong and Zhang [38]. The tests of the first two papers are based on a (weighted) residual sum of squares and are in the spirit of Härdle and Mammen [17]. The asymptotic normality of the test statistic is achieved by a clever application of central limit theorems for sequences of degenerate $U$-statistics. More precisely, Fan and Li [14] (FL) based their test on a quadratic form of the estimated residuals. Since no rank transformation is involved, they had to weight each residual with estimators of marginal and high-dimensional densities, to get rid of the denominator in the Nadaraya–Watson estimator. Consequently, two different smoothing parameters need to be involved. It is heuristically argued that local alternatives only can be detected when they approach the null model at the rate $O((nh^{p/2})^{-1/2})$, which gets worse as the dimension of $X$ increases. The estimator of $\beta$, being square-root consistent, does not have any impact on the limit distribution because the other quantities converge at a slower rate, thus compensating for the effect of estimating unknown parameters. In a general situation of testing a model or hypothesis, efficient methods involve test statistics and estimators which admit expansions of the same order. See, for example, [9], to name only one landmark paper on this topic. Unless some orthogonality assumptions are satisfied, the parameter estimator does have an impact on the limit, and martingale transformations, as in [36], were designed to keep track of this issue. See also [34]. Efficient model checks would therefore create terms which when replacing $\hat{\beta}$ with $\beta$ are not negligible and thus have an impact on the distributional behavior of the test statistic. As to practical applications, computation of critical values would then not be easy. Worse than that, the complicated geometric structure of the test statistic would not enable us to derive optimal scores. Actually, these are only two of several reasons why we designed our test as we did. There are others. As a by-product, the assumptions on the design variable $X$ remain weak. No additional support or higher smoothness conditions need to be assumed. The variable $Y$ may be discrete and no joint density of $X$ and $Y$ is required. Compared with Fan and Li [14], Aït-Sahalia, Bickel and Stoker [1] is mainly concerned with the problem of dimension reduction for high-dimensional inputs. Only some comments on the applicability to single-index models are included. Their test statistic is a sum of weighted residual squares, the weights now being deterministic functions of the regressors. In their Proposition 2 the local power of the test is derived when the alternatives tend to the null model at a rate depending on $p$. It should also be mentioned that the test statistic admits a bias increasing to infinity as $n \to \infty$. Moreover, the constants defining the asymptotic bias are unknown and require further smoothing when being estimated. Similarly, in Xia, Li, Tong and Zhang [38], who extended the marked empirical process approach of Stute, González Manteiga and Presedo Quindimil [35] in the parametric case to the single index model. Compared with these papers our test achieves local power known from parametric tests, though the nonparametric components can only be estimated at a worse rate. Mathematically, we have to pay a price for this. For example, Theorem 2.1 cannot be obtained by just applying Taylor's expansion and $U$-statistic theory. Rather, our proofs require some new techniques involving (local and global) properties of the rank-transformed projected values $\hat{\beta}^T X_i, 1 \leq i \leq n$. Unfortunately, techniques also developed in [31] to analyze the (rank-transformed) nearest-neighbor regression function estimator at a point are of no help here.

## 3. Simulation study and applications.

3.1. *A simulation study.* In our simulations we studied two models. The first is with continuous response, namely,

$$(3.1) \qquad Y = (\beta^T X)^3 + c\left(\sum_{l=1}^{p} |x_l|\right) + \epsilon,$$

where $X$ and $\epsilon$ are independent, $x_l$ are the components of $X$ and the distributions of $X$ and $\epsilon$ are $N(0, I_p)$ and $N(0, 1)$, respectively. The hypothetical model is $\Phi(\beta^T X) = (\beta^T X)^3$ and $s(X) = \sum_{l=1}^{p} |x_l|$. Therefore, the null model holds if and only if $c = 0$.

The second model is with binary response,

$$
\begin{aligned}
(3.2) \quad Y &= \frac{\exp(-\beta^T X + c(\sum_{l=1}^{p} |x_l|))}{1 + \exp(-\beta^T X + c(\sum_{l=1}^{p} |x_l|))} + \epsilon \\
&=: \Phi(\beta^T X + cs(X)) + \epsilon,
\end{aligned}
$$

where $Y = 0, 1$ is a binary variable for which $Y = 1$ with probability $\Phi(\beta^T \mathbf{x} + cs(\mathbf{x}))$ for any given $X = \mathbf{x}$. Also, here $c = 0$ corresponds to the hypothetical model, that is, the logit model. It is heteroscedastic, and $X$ and $\epsilon$ are not independent. Again, $X \sim \mathcal{N}(0, I_p)$. We used $c = 1, 2, 3$ to investigate the power of the test.

Two weight functions were considered in the simulation, $W_1(\mathbf{x}) = s(\mathbf{x})$ and $W_2(\mathbf{x}) = \sum_{l=1}^{p} x_l^2$. Based on our findings in Section 2, $W_1$ is optimal for model (3.1) as $\epsilon$ is independent of $X$, and $W_2$ is a natural candidate for an even function. For model (3.2), we also use these two weight functions due to the following observation: When $c$ is small, $\Phi(-\beta^T \mathbf{x} + cs(\mathbf{x}))$ is close to $\Phi(-\beta^T \mathbf{x}) + c\Phi'(\beta^T \mathbf{x})s(\mathbf{x})$, where $\Phi'(\cdot)$ is the derivative of $\Phi(\cdot)$. Therefore, $s(\mathbf{x})$ is also a good choice of a weight function in this case.

In order to implement the omnibus test based on $\tilde{T}_n = \sup_\gamma |\hat{T}_n(\gamma)|$ of Theorem 2.3, we have to use a resampling approximation to determine critical values. The wild bootstrap is clearly an option. In view of (2.6), however, we suggest the following algorithm: for any $\gamma$, $\hat{T}_n(\gamma)$ is asymptotically equal to $\mu + n^{-1/2} \sum_{i=1}^{n} \varepsilon_i [W_i - \bar{W}(U_i)]$. Under $H_0$, $\mu = 0$. For any i.i.d. random variables $e_i, i = 1, \ldots, n$, independent of the $(x_i, y_i)$'s with mean 0 and variance 1, it is easy to prove that, for almost all sequences $\{(x_1, y_1), \ldots, (x_n, y_n), \ldots\}$, the process $T_n^r(\gamma) = n^{-1/2} \sum_{i=1}^{n} e_i \hat{\varepsilon}_i [W_i - \bar{W}_n^{(i)}(\hat{U}_i)]$ has the same limit as $\hat{T}_n(\gamma)$. It is worthwhile noting that, using this resampling scheme, we do not need to estimate the variance. In a different setup, this algorithm has been used by Zhu [40] and Zhu and Ng [42]. The proof and the procedure are similar. We omit the details. To implement the test, we can generate, by Monte Carlo, $m$ sets of $\{e_1, \ldots, e_n\}$ and then compute $m$ values of $\tilde{T}_n^r = \sup_\gamma |\hat{T}_n^r(\gamma)|$. The $[(1 - \alpha)m]$th value can be used as the critical value, where $\alpha$ is the significance level and $[a]$ stands for the integer part of $a$. In the following simulation, we used standard normal random variables $e_i$.

Another concern is bandwidth selection. As we noticed in Remark 2.5, $h \sim n^{-1/3}$. In other words, compared with nonparametric estimation of regression, in the context of model checking, undersmoothing is needed. So existing bandwidth selection methods cannot be recommended in the setting of this paper and, indeed, may lead to a considerable bias. Therefore, we adopt a semidata driven selection procedure. The steps are as follows:

TABLE 1
*Size of the tests $\bar{T}_n$ and $\tilde{T}_n$*[a]

|  |  | Model (3.1) | | |  | Model (3.2) | |
| --- | --- | --- | --- | --- | --- | --- | --- |
|  |  | $n=50$ | $n=100$ |  |  | $n=50$ | $n=100$ |
| $W_1$ | $p=2$ | 0.048(0.046) | 0.045(0.047) | $W_1$ | $p=2$ | 0.060(0.056) | 0.057(0.054) |
| $W_1$ | $p=3$ | 0.053(0.053) | 0.047(0.052) | $W_1$ | $p=3$ | 0.054(0.052) | 0.052(0.054) |
| $W_2$ | $p=2$ | 0.048(0.047) | 0.052(0.053) | $W_2$ | $p=2$ | 0.055(0.054) | 0.052(0.054) |
| $W_2$ | $p=3$ | 0.047(0.053) | 0.046(0.051) | $W_2$ | $p=3$ | 0.055(0.055) | 0.050(0.054) |
| $\tilde{T}_n$ | $p=2$ | 0.048(0.051) | 0.053(0.052) | $\tilde{T}_n$ | $p=2$ | 0.058(0.056) | 0.057(0.051) |
| $\tilde{T}_n$ | $p=3$ | 0.045(0.048) | 0.052(0.054) | $\tilde{T}_n$ | $p=3$ | 0.061(0.053) | 0.054(0.049) |

[a] The values in parentheses are the estimated sizes when the bandwidth is selected by a grid search.

1. Select $h_1$ by minimizing the mean integrated squared error, subject to weight function $W(\cdot)$,

$$(3.3) \qquad MISE(h) = \sum_{j=1}^{n}(Y_j - \hat{\psi}_n^{(j)}(\hat{U}_j))^2 W(X_j)^2,$$

which is analogous to the criterion used by Härdle, Hall and Ichimura [16]. The kernel $K$ is $15/16(1-u^2)^2 I(|u| \leq 1)$; see [17].

2. Our final choice for $h$ is $h = h_1 \times n^{-1/3+1/5}$.

The rationale of this algorithm is that, under our conditions and the choice of the kernel function, the rate of $h_1$ is $n^{-1/5}$. Therefore, $h$ is of the order $n^{-1/3}$ and, hence, ensures convergence of the test statistic. For validation purposes we also considered a grid point search and chose $h$ so that the empirical level was closest to the nominal level.

Finally, we need to estimate the parameter $\beta$. There are at least three methods in the literature; see [16, 20, 25]. In our simulation study we applied Li and Duan's least squares estimator for ease of implementation.

We considered the case with $p=2,3$ and $\beta = (1,-1)^T/\sqrt{2}$, $\beta = (1,-1,1)^T/\sqrt{3}$, respectively. The sample sizes were $n = 50, 100$. The significance level was $\alpha = 0.05$. The test statistics were computed for 1000 replications.

Table 1 presents the attained levels for the various scenarios.

It becomes apparent that the significance level is well attained in most cases, although, for model (3.2), the size of the tests for $n = 50$ is slightly larger than 0.05. Furthermore, the size of the tests with the bandwidth selected by the above algorithm is similar to that obtained from the grid point search. This shows that our data-driven approach works well. We will therefore use this algorithm also to select the bandwidth in the following simulation and the applications to two real data examples.

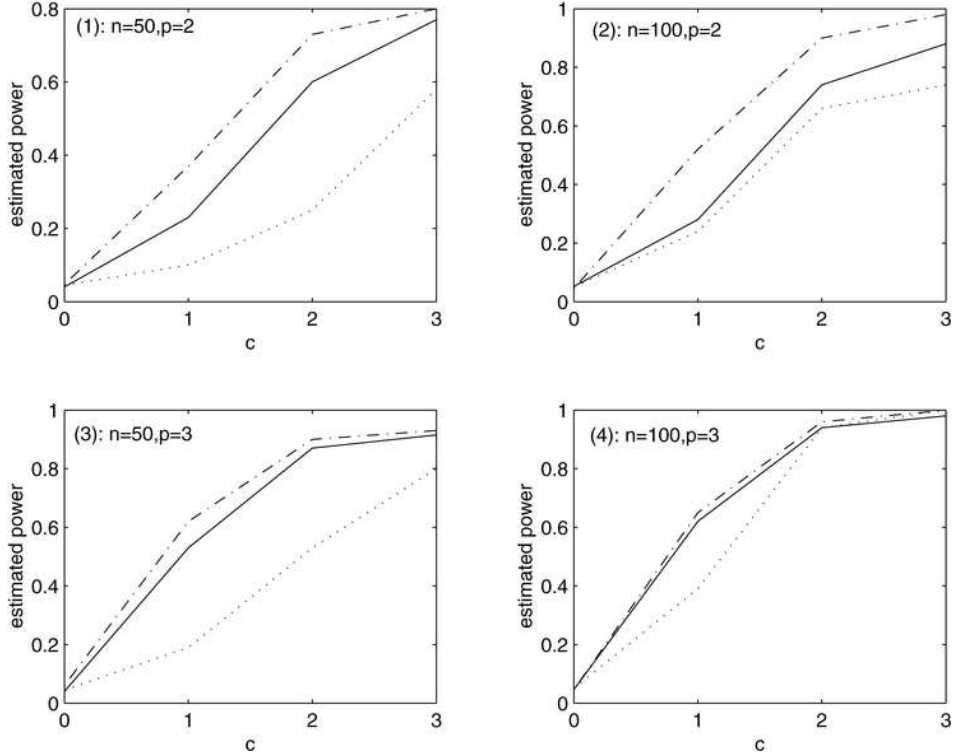

FIG. 1. *The estimated power for model* (3.1): *The dashdot line is for our test with the weight function* $W_1$, *the solid line with the weight function* $W_2$, *and the dotted line is for* $\tilde{T}_n$.

To demonstrate power through simulations, we considered models (3.1) and (3.2) with $c = 1, 2, 3$.

For model (3.1), as expected, the test $T_n$ based on the optimal $W_1$ outperforms the others. In model (3.2), when we have dependent errors and $T_n$ is no longer optimal, all three tests have a similar behavior.

To compare the performance of our method with other existing tests through a simulation study, we considered two scenarios. The first aim was to test the single index model versus the existence of interaction effects. Particularly, we considered

$$(3.4) \qquad m(\mathbf{x}) = (\beta^T \mathbf{x})^3 + c_1|x_1 x_2| + c_2|x_1 x_3| + c_3|x_2 x_3|.$$

For nonvanishing $c$'s, this model allows for interaction terms. The comparison is among our maximin test, the omnibus test $\tilde{T}_n$, Fan and Li [14] (FL-test) and Aït-Sahalia, Bickel and Stoker [1] (ABS-test). In the simulation, similar to the previous case, we took $\beta = (1, -1, 1)^T / \sqrt{3}$. The sample size was $n = 50$, while the significance level was 0.05. The constants were taken to be equal: $c_1 = c_2 = c_3 = c$ with $c = 0, 1.0, 2.0, 3.0$. $c \neq 0$ corresponds to the alternative. In Figure 3 the estimated power was computed from 1000 replications. Recall that FL- and ABS-tests require selection of two bandwidths. Since the significance levels of their tests heavily depend on the choice of the bandwidths and there is no data driven selection, a fair comparison causes some problems. In a simulation study, however, one may determine (through replications) the bandwidth on a grid in such a way that the nominal level is best attained. In this way we are able to produce tests which attain the right level for the null model.

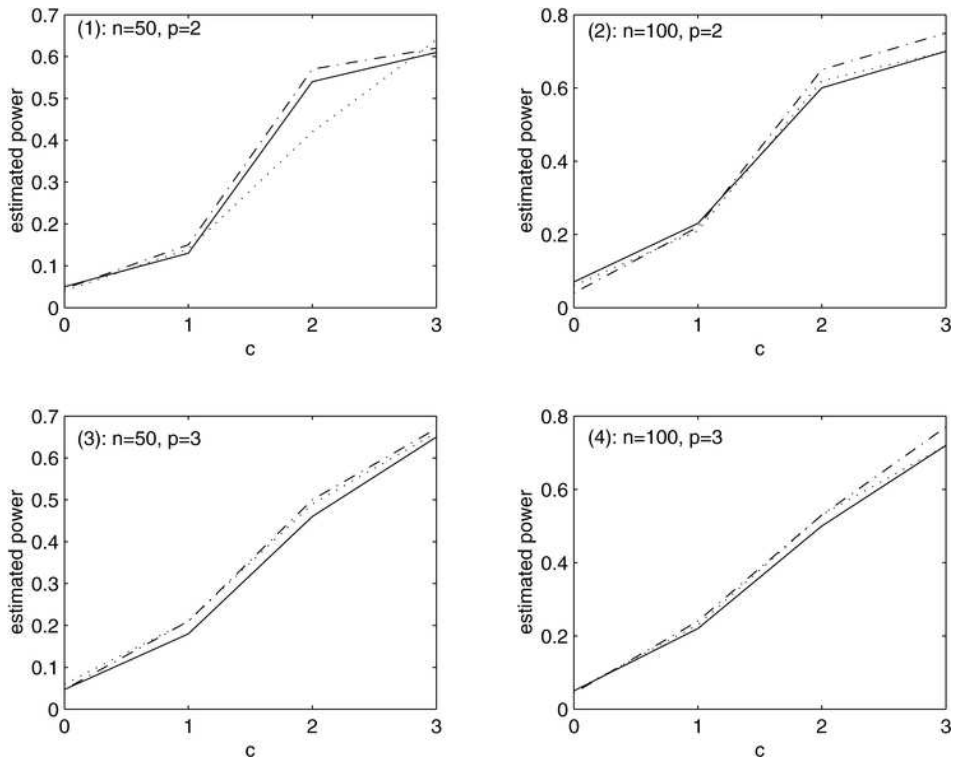

FIG. 2. *The estimated power for model* (3.2): *The dashdot line is for our test with the weight function $W_1$, the solid line with the weight function $W_2$, and the dotted line is for $\tilde{T}_n$.*

We also ran many simulations with other bandwidths. It turned out that the FL-test and the ABS-test are nonrobust in $h$ so that the nominal level may not be attained after a slight change in $h$.

As expected, $T_n$ with optimal weight $W_1$ has larger power than the test with weight function $W_2$. $\tilde{T}_n$ has a power similar to $T_n$ with $W_2$. The FL- and ABS-tests are clearly outperformed but behave similarly otherwise in the situation considered by us. Similar to the case with model (3.1), the FL-test has larger power than the ABS-test.

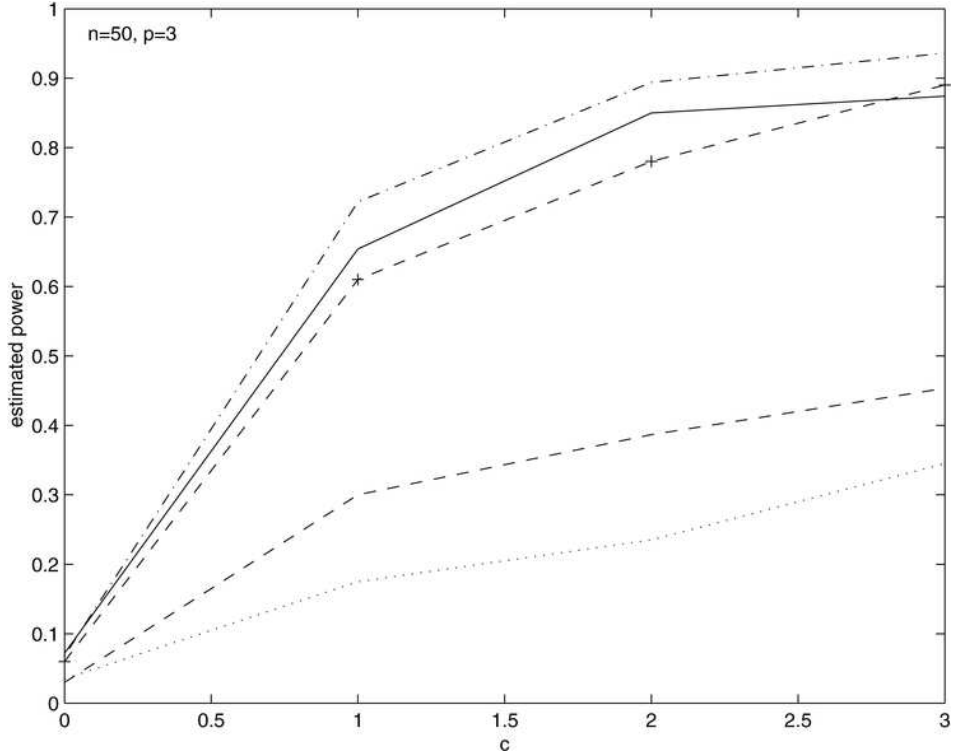

FIG. 3. *The estimated power for model* (3.4): *The dashdot line is for the maximin test with the weight function* $W_1$, *the solid line with the weight function* $W_2$; *the dotted line is for the ABS-test, the dashed line for the FL-test, and the dashed line plus star* $*$ *for* $\tilde{T}_n$.

We also compared the performance of all tests for a model studied by Xia, Li, Tong and Zhang [38] in their Example 1, where, in our notation, $p = 2$ and

$$m(x) = x_1 + x_2 + 4\exp\{-(x_1 + x_2)^2\} + c(x_1^2 + x_2^2)^{1/2},$$

and the errors $\varepsilon$ are independent of $X$ with $\varepsilon \sim \mathcal{N}(0, \sigma_\varepsilon^2)$.

In Table 2 we report on the power results of $T_n$ with $W_1(\cdot)$ and $W_2(\cdot)$, $\tilde{T}_n$, ABS- and FL-tests and the XLTZ-test. The bootstrap approximation of the XLTZ-test is similar to that of Theorem 2.3. For $T_n$, we again used the weights $W_1(x) = |x_1| + |x_2|$ and $W_2(x) = x_1^2 + x_2^2$. The significance level was 0.05. The test statistics were computed for 1000 replications. Note that these two weights are not optimal for this model. We do not report the results with the optimal weights because the previous simulations have provided evidence of its good performance and, from Table 2, we can see that the suboptimal weights $W_1$ and $W_2$ already work well. Again, for ABS and FL, bandwidths were chosen so as to yield the nominal level under $H_0$ as closely as possible.

TABLE 2
*Estimated power of six tests with $n = 50, p = 2$*[a]

| $\sigma_\varepsilon$ | 0.30 | | | 0.50 | | |
|---|---|---|---|---|---|---|
| c | 0 | 0.25 | 0.50 | 0 | 0.25 | 0.50 |
| $T_n(W_1)$ | 0.044 | 0.122 | 0.508 | 0.052 | 0.106 | 0.452 |
| $T_n(W_2)$ | 0.060 | 0.092 | 0.408 | 0.062 | 0.090 | 0.300 |
| XLTZ-test | 0.063 | 0.099 | 0.376 | 0.043 | 0.043 | 0.163 |
| $\tilde{T}_n$ | 0.063 | 0.090 | 0.350 | 0.043 | 0.073 | 0.253 |
| ABS-test | 0.050 | 0.060 | 0.140 | 0.050 | 0.055 | 0.085 |
| FL-test | 0.042 | 0.052 | 0.090 | 0.050 | 0.046 | 0.065 |

[a] $T_n(W_i)$, $i = 1, 2$, stand for the tests $T_n$ with $W_1$ and $W_2$, respectively.

In Table 2, the values for the XLTZ-test are from Table 1 of [38]. We see that $T_n$ with $W_1$ is best. Second, between $\tilde{T}_n$ and the XLTZ-test, when the variance $\sigma_\varepsilon^2$ of the errors $\varepsilon_i$ is small, the XLTZ-test is slightly better, while when $\sigma_\varepsilon^2$ gets large, $\tilde{T}_n$ outperforms the XLTZ-test. Third, comparing $\tilde{T}_n$ with $T_n$ with $W_2$, we see that $\tilde{T}_n$ performs slightly worse. For this model, we find that the ABS- and FL-tests do not work well.

3.2. *Applications.* In this section we apply our test to two data sets.

EXAMPLE 3.1. The data set is the bull data; see [24]. The data are the measured characteristics of 76 young bulls sold at an auction. It is interesting to study the relationship between the selling prices and the characteristics of the bulls: yearling height at shoulder; fat-free body (pounds); percentage of fat-free body; scale from 1 (small) to 8 (large); back fat (inches); sale height at shoulder (inches) and scale weight (pounds). The response $Y$ is the standardized selling price and the other standardized measurements are the covariates $X = (x_1, \ldots, x_7)$. Figure 4(a) provides a plot of $\hat{\beta}^T X$ against the response $Y$. This linear fitting was also used in [24]. There is some indication of a relationship between the residuals $\hat{\epsilon}_j$ and $\hat{\beta}^T X_j$, see Figure 4(b). We tested the linearity of the model using the Stute, González Manteiga and Presedo Quindimil [35] test. The $p$-value was 0.044. Therefore, the linear model needs to be rejected at level $\alpha = 0.05$.

Next consider single-index fitting. Again $\beta$ was estimated as in [25]. To justify their estimation method, we first tested the elliptical symmetry of the distribution of $X$. The nonparametric Monte Carlo test proposed by Zhu and Neuhaus [41] was employed. The $p$-value was 0.83. The statistic $\bar{T}_n$ was computed for the weight function $W(\mathbf{x}) = \sum_{j=1}^p x_j^2$. The kernel function

SINGLE-INDEX MODELS 21

$K(\cdot)$ is the same as for (3.3), and the bandwidth is $h = 0.35$. The $p$-value was 0.310. Therefore, a single-index model need not be rejected.

EXAMPLE 3.2. The data are the automobile collision data as analyzed by Härdle, Hall and Ichimura [16]. The sample size is $n = 58$. We also tested the elliptical symmetry of the distribution of the $X$-data using the nonparametric Monte Carlo test of Zhu and Neuhaus [41]. The $p$-value was 0.25. This justifies the use of the Li–Duan method for estimating the projection direction $\beta$ for this data set. For a single-index fitting, the kernel function $K(\cdot)$ was again the same as for (3.3), the bandwidth was $h = 0.4$, while the weight function was $W(\mathbf{x}) = \sum_{j=1}^{p} x_j^2$. The test statistic $\bar{T}_n$ was used and the asymptotic $p$-value was 0.32. The single-index model is therefore tenable.

**4. Proofs.** To prove Theorem 2.1, we expand our test statistic $\hat{T}_n$ as

$$
\begin{aligned}
n^{1/2}\hat{T}_n &= \sum_{j=1}^{n} \hat{\varepsilon}_j W_j = \sum_{j=1}^{n} [Y_j - \psi_n^{(j)}(F_n(\hat{\beta}^T X_j))] W_j \\
(4.1) \quad &= \sum_{j=1}^{n} [Y_j - Y_j^0 - \psi_n^{(j)}(F_n(\hat{\beta}^T X_j)) + \psi_{n0}^{(j)}(F_n(\hat{\beta}^T X_j))] W_j \\
&\quad + \sum_{j=1}^{n} [Y_j^0 - \psi_{n0}^{(j)}(F_n(\hat{\beta}^T X_j))] W_j \equiv I + II,
\end{aligned}
$$

where $Y_j^0$ is computed under the null model $s \equiv 0$, and $\psi_{n0}^{(j)}$ is computed as $\psi_n^{(j)}$, with the same $\hat{\beta}$ but with $Y_j^0$. The second sum will be further

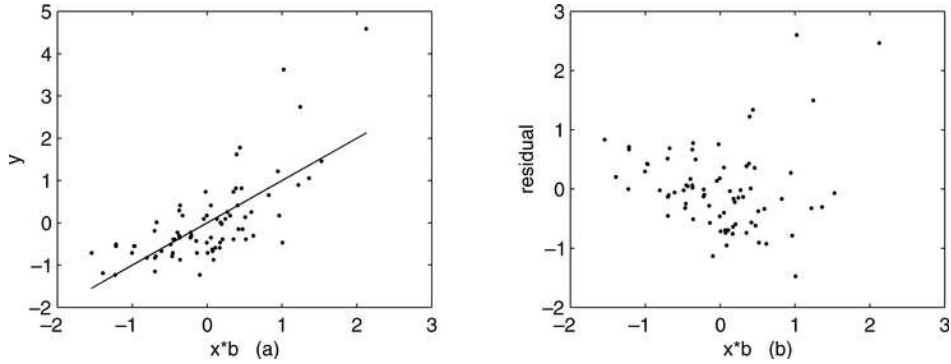

FIG. 4. (a) *Fit to the bulls data: the projected data $\hat{\beta}^T X_j$ versus the linear fit (solid line) and the response data (dots);* (b) *the projected data versus the residuals.*

decomposed. For this, put

$$\bar{\psi}_{n0}^{(j)}(u) = \frac{1}{(n-1)h} \sum_{\substack{i=1 \\ i \neq j}}^{n} Y_i^0 K\left(\frac{u - F(\beta^T X_i)}{h}\right).$$

This function is based on the true $\beta$ and $F$ and is therefore unknown in practice. It will, however, play an important role in proofs, since it is close to $\psi_{n0}^{(j)}$ and, on the other hand, is computed from independent observations. Write

$$II = \sum_{j=1}^{n} [Y_j^0 - \bar{\psi}_{n0}^{(j)}(F(\beta^T X_j))] W_j$$

$$+ \sum_{j=1}^{n} [\bar{\psi}_{n0}^{(j)}(F(\beta^T X_j)) - \psi_{n0}^{(j)}(F_n(\hat{\beta}^T X_j))] W_j \equiv III + IV.$$

Observe that

$$III = \sum_{j=1}^{n} Y_j^0 W_j - \frac{1}{(n-1)h} \sum_{\substack{j=1 \\ i \neq j}}^{n} \sum_{i=1}^{n} Y_i^0 W_j K\left(\frac{U_j - U_i}{h}\right),$$

with

$$U_j = F(\beta^T X_j), \qquad j = 1, \ldots, n,$$

being independent and uniformly distributed on $[0,1]$. Hence, $III$ is a $U$-statistic of degree two. Summarizing, we have

(4.2) $$\sum_{j=1}^{n} \hat{\varepsilon}_j W_j = I + III + IV.$$

After standardization, term $I$ will be shown to tend to a limit which depends on the shift $s$ and, hence, will determine the local power of the test. As already mentioned, $III$ is a $U$-statistic of degree two, with a kernel depending on $h$, and hence on $n$. The term $IV$ is more complicated, since the kernel contains empirical quantities. After all, it will turn out that $III$ and $IV$ admit i.i.d. representations which will partly cancel out and jointly determine the (limit) distribution of $\hat{T}_n$ under $H_0$. To carry out this program, note that both $\psi_n^{(j)}$ and $\psi_{n0}^{(j)}$ are evaluated at $F_n(\hat{\beta}^T X_j)$. Hence, the mathematical analysis of our test statistic requires a careful study of the terms

(4.3) $$K\left(\frac{F_n(\hat{\beta}^T X_j) - F_n(\hat{\beta}^T X_i)}{h}\right), \qquad 1 \leq i \neq j \leq n.$$

For this, denote by $F_n^\theta$ the empirical distribution function of $\theta^T X_1, \ldots, \theta^T X_n$. Hence, $F_n = F_n^\theta$ if $\theta = \hat{\beta}$. Since $K$ has compact support, say $[-1, 1]$, indices $i, j$ only contribute to (4.3) if

$$\text{(4.4)} \qquad |F_n^\theta(\theta^T X_j) - F_n^\theta(\theta^T X_i)| \leq h, \qquad \theta = \hat{\beta}.$$

Since by assumption B(iv)

$$n^{1/2}(\hat{\beta} - \beta) = O_\mathbb{P}(1),$$

for each given $\varepsilon > 0$, we may find a large constant $C$ such that

$$\mathbb{P}(n^{1/2}\|\hat{\beta} - \beta\| \geq C) \leq \varepsilon \qquad \text{for all } n \geq 1.$$

In other words, up to a small event, $\hat{\beta}$ is contained in the $Cn^{-1/2}$-neighborhood of $\beta$. The first goal will be to analyze the effect of replacing $F_n(\hat{\beta}^T X_j)$ and $F_n(\hat{\beta}^T X_i)$ in (4.3) with $U_j = F(\beta^T X_j)$ and $U_i = F(\beta^T X_i)$, respectively, subject to (4.4). Introduce $F^\theta$, the distribution function of $\theta^T X$. Hence, $F = F^\theta$ for $\theta = \beta$.

In our first lemma we derive a maximal bound for $F^\theta - F^\beta$ evaluated at $\theta^T X_j$ and $\beta^T X_j$. Recall that, by assumption B(i), $\mathbb{E}\|X\|^\gamma < \infty$. This implies that

$$\max_{1 \leq i \leq n} \|X_i\| = O_\mathbb{P}(n^\alpha) \qquad \text{for } \alpha = \gamma^{-1}.$$

For this reason, it will suffice to analyze all leading and error terms on the set where

$$\text{(4.5)} \qquad \max_{1 \leq i \leq n} \|X_i\| \leq C_1 n^\alpha \qquad \text{for some large finite } C_1.$$

Denote by $\Theta$ the set of all $p$-vectors.

LEMMA 4.1. *Put, for each $\theta \in \Theta$ and $1 \leq j \leq n$,*

$$a_j^\theta := F^\theta(\theta^T X_j) - F^\beta(\beta^T X_j).$$

*We then have, on the set* (4.5),

$$\max_{\|\theta - \beta\| \leq Cn^{-1/2}} \max_{1 \leq j \leq n} |a_j^\theta| = O_\mathbb{P}(n^{-1/2+\alpha}).$$

PROOF. We shall first deal with an upper bound for the $a_j^\theta$'s. Fix a possible value $\mathbf{x}_j$ of $X_j$. Then

$$a_j^\theta = F^\theta(\theta^T \mathbf{x}_j) - F^\beta(\beta^T \mathbf{x}_j) = \mathbb{P}(\theta^T X \leq \theta^T \mathbf{x}_j) - \mathbb{P}(\beta^T X \leq \beta^T \mathbf{x}_j)$$

$$= \mathbb{P}(\theta^T X \leq \theta^T \mathbf{x}_j, \beta^T X \leq \beta^T \mathbf{x}_j) + \mathbb{P}(\theta^T X \leq \theta^T \mathbf{x}_j, \beta^T X > \beta^T \mathbf{x}_j)$$

$$- \mathbb{P}(\beta^T X \leq \beta^T \mathbf{x}_j) \leq \mathbb{P}(\theta^T X \leq \theta^T \mathbf{x}_j, \beta^T X > \beta^T \mathbf{x}_j).$$

Now, $\theta^T X \leq \theta^T \mathbf{x}_j$ implies

$$\beta^T X = \theta^T X + (\beta - \theta)^T X \leq \theta^t \mathbf{x}_j + (\beta - \theta)^T X$$
$$= \beta^T \mathbf{x}_j + (\beta - \theta)^T (X - \mathbf{x}_j) \leq \beta^T \mathbf{x}_j + Cn^{-1/2}\{\|X\| + \|\mathbf{x}_j\|\}.$$

Under (4.5) we therefore obtain, for each $1 \leq j \leq n$,

$$a_j^\theta \leq \mathbb{P}(\beta^T \mathbf{x}_j < \beta^T X \leq \beta^T \mathbf{x}_j + 2CC_1 n^{-1/2+\alpha}) + \mathbb{P}(\|X\| > C_1 n^\alpha).$$

Since, by B(ii), $\beta^T X$ has a bounded density, the first probability is $O(n^{-1/2+\alpha})$. As to the second probability, apply Markov's inequality to get

$$\mathbb{P}(\|X\| > C_1 n^\alpha) \leq \frac{\mathbb{E}\|X\|^\gamma}{C_1^\gamma n}.$$

This completes the proof. For the lower bound, just reverse the roles of $\theta$ and $\beta$. Now one needs the fact that the densities of $\theta^T X$ are uniformly bounded for all $\theta$ in a small neighborhood of $\beta$. □

In the following lemma we investigate the local oscillations of the empirical process

$$(x, \theta) \to F_n^\theta(x) - F^\theta(x)$$

in a neighborhood of $\beta$. For this, introduce

$$G_n^\theta(x, y) := F_n^\theta(x) - F^\theta(x) - F_n^\beta(y) + F^\beta(y)$$

for $\theta \in \Theta$ and $x, y \in \mathbb{R}$ satisfying

(i) $\|\theta - \beta\| \leq Cn^{-1/2}$,
(ii) $|x - y| \leq C_1 n^{-1/2+\alpha}$.

LEMMA 4.2. *Under the assumptions of Theorem* 2.1, *we have*

$$\sup_{x,y;\theta} |G_n^\theta(x, y)| = O_\mathbb{P}(\sqrt{n^{-3/2+\alpha} \ln n}),$$

*where the supremum extends over all $x, y$ and $\theta$ satisfying* (i) *and* (ii).

PROOF. The proof is a modification of the proof of Theorem 37 in [26], page 34. First note that the halfspaces form a class with a polynomial covering number. The measure of each set involved in the above supremum, $F^\theta(x) - F^\beta(y)$, is bounded from above in absolute value by

$$|\mathbb{P}(\theta^T X \leq x) - \mathbb{P}(\beta^T X \leq y)|$$
$$\leq |\mathbb{P}(\theta^T X \leq x) - \mathbb{P}(\beta^T X \leq x)|$$
$$+ |\mathbb{P}(\beta^T X \leq x) - \mathbb{P}(\beta^T X \leq y)| \leq C_2 n^{-1/2+\alpha},$$

by (i), (ii) and assumption B. For the first difference apply a technique already used in the proof of the previous lemma. If we replace the small $\varepsilon$ in Pollard's [26] Theorem 37 by a large $K > 0$ and set

$$\alpha_n^2 = \frac{\ln n}{n\delta_n^2}$$

therein, we obtain the required in-probability bound $O(\delta_n^2 \alpha_n)$, rather than a convergence rate to zero. Here $\delta_n^2$ equals the maximal measure of the included sets. Since $\delta_n^2 = O(n^{-1/2+\alpha})$, the result follows. $\square$

In the next lemma, we expand $n^{-1/2} III$ into a sum of independent random variables plus a negligible error. The leading term will contribute to the limit of our test statistic when the null hypothesis is true. Recall

$$\bar{W}(u) = \mathbb{E}[W_1 | U_1 = u].$$

LEMMA 4.3. *Under the assumptions of Theorem* 2.1, *we have in probability as* $n \to \infty$,

$$n^{-1/2} III \equiv S_{n_1} = n^{-1/2} \sum_{j=1}^{n} \varepsilon_j W_j - n^{-1/2} \sum_{j=1}^{n} [Y_j^0 \bar{W}(U_j) - \mathbb{E} Y_1^0 \bar{W}(U_1)] + o_{\mathbb{P}}(1)$$

$$= n^{-1/2} \sum_{j=1}^{n} \{\varepsilon_j [W_j - \bar{W}(U_j)]$$

$$- \Phi(\beta^T X_j) \bar{W}(U_j) + \mathbb{E}[\Phi(\beta^T X_j) \bar{W}(U_j)]\}$$

$$+ o_{\mathbb{P}}(1).$$

PROOF. $S_{n_1}$ is a $U$-statistic of degree two with a kernel depending on $h$ and therefore on $n$. The Hájek projection of $Y_i^0 W_j K(\frac{U_j - U_i}{h})$ equals

$$Y_i^0 \int_0^1 \bar{W}(v) K\left(\frac{v - U_i}{h}\right) dv + W_j \int_0^1 \psi(u) K\left(\frac{U_j - u}{h}\right) du$$

$$- \int_0^1 \int_0^1 \bar{W}(v) \psi(u) K\left(\frac{v - u}{h}\right) dv\, du.$$

Conclude that the Hájek projection of $S_{n_1}$ equals

$$\hat{S}_{n_1} = n^{-1/2} \sum_{j=1}^{n} Y_j^0 W_j - n^{-1/2} h^{-1} \sum_{i=1}^{n} Y_i^0 \int_0^1 \bar{W}(v) K\left(\frac{v - U_i}{h}\right) dv$$

$$- n^{-1/2} h^{-1} \sum_{j=1}^{n} W_j \int_0^1 \psi(u) K\left(\frac{U_j - u}{h}\right) du$$

$$+ h^{-1} n^{1/2} \int_0^1 \int_0^1 \bar{W}(v) \psi(u) K\left(\frac{v - u}{h}\right) dv\, du.$$

Furthermore (see [27]),
$$\mathbb{E}\{S_{n_1} - \hat{S}_{n_1}\}^2 = O\left(\frac{1}{nh}\right),$$
whence
$$S_{n_1} - \hat{S}_{n_1} = O_{\mathbb{P}}((nh)^{-1/2}) = o_{\mathbb{P}}(1).$$

Hence, it suffices to further expand $\hat{S}_{n_1}$. For this, put
$$E_h = \int_0^1 \int_0^1 \bar{W}(u)\psi(v)K\left(\frac{v-u}{h}\right)dv\,du$$
and consider
$$R_{n_1} = n^{-1/2}h^{-1}\sum_{i=1}^n \left[Y_i^0 \int_0^1 \bar{W}(v)K\left(\frac{v-U_i}{h}\right)dv - E_h\right]$$
$$+ n^{-1/2}h^{-1}\sum_{j=1}^n \left[W_j \int_0^1 \psi(u)K\left(\frac{U_j-u}{h}\right)du - E_h\right]$$
$$- n^{-1/2}\sum_{i=1}^n [Y_i^0 \bar{W}(U_i) - \mathbb{E}(Y_1^0 \bar{W}(U_1))]$$
$$- n^{-1/2}\sum_{j=1}^n [W_j \psi(U_j) - \mathbb{E}(W_1 \psi(U_1))].$$

It may be written as a single sum of centered i.i.d. random variables. Its variance is bounded from above by the second moment of
$$Y_1^0 \left[h^{-1}\int_0^1 \bar{W}(v)K\left(\frac{v-U_1}{h}\right)dv - \bar{W}(U_1)\right]$$
$$+ W_1 \left[h^{-1}\int_0^1 \psi(u)K\left(\frac{U_1-u}{h}\right)du - \psi(U_1)\right],$$
which is easily seen to go to zero as $h \to 0$. Conclude that $R_{n_1} = o_{\mathbb{P}}(1)$ and, therefore,
$$S_{n_1} = \hat{S}_{n_1} + o_{\mathbb{P}}(1)$$
$$= n^{-1/2}\sum_{j=1}^n Y_j^0 W_j - n^{-1/2}\sum_{j=1}^n [Y_j^0 \bar{W}(U_j) - \mathbb{E}(Y_1^0 \bar{W}(U_1))]$$
$$- n^{-1/2}\sum_{j=1}^n [W_j \psi(U_j) - \mathbb{E}(W_1 \psi(U_1))] - n^{1/2}h^{-1}E_h + o_{\mathbb{P}}(1)$$
$$= n^{-1/2}\sum_{j=1}^n \varepsilon_j W_j - n^{-1/2}\sum_{j=1}^n [Y_j^0 \bar{W}(U_j) - \mathbb{E}(Y_1^0 \bar{W}(U_1))]$$

$$+ n^{1/2}[\mathbb{E}(W_1\psi(U_1)) - h^{-1}E_h] + o_{\mathbb{P}}(1).$$

To complete the proof of the lemma, it suffices to show, in view of assumption C(i), that the last bracket is $O(h^2)$. But

$$[\cdots] = \int_0^1 \bar{W}(v)\left[\psi(v) - h^{-1}\int_0^1 \psi(u)K\left(\frac{v-u}{h}\right)du\right]dv$$

$$= \int_0^1 \bar{W}(v)\left[\psi(v) - \int_{(v-1)/h}^{v/h} \psi(v-sh)K(s)\,ds\right]dv.$$

For $h \leq v \leq 1-h$, the inner integral extends over the whole support of $K$, namely, $[-1,1]$. Using the facts that $K$ is symmetric at zero, $\int_{-1}^1 K(s)\,ds = 1$ and $\psi$ is twice continuously differentiable, Taylor's expansion yields that the difference is uniformly in $h \leq v \leq 1-h$ of the order $O(h^2)$. For $0 \leq v < h$ (and similarly for $1-h < v \leq 1$), the difference is $O(h)$. Since, however, $0 \leq v < h$ has Lebesgue measure $h$, we also obtain the upper bound $h^2$ for this part of the integral. $\square$

The quantity $S_{n_2}$ introduced and studied below will be the leading term for $n^{-1/2}I$ with $I$ from the expansion (4.2).

LEMMA 4.4. *Under the assumptions of Theorem 2.1, we have in probability as $n \to \infty$,*

$$S_{n_2} \equiv n^{-1}\sum_{j=1}^n s(X_j)W_j - \frac{1}{n(n-1)h}\sum_{j=1}^n \sum_{\substack{i=1 \\ i \neq j}}^n s(X_i)W_j K\left(\frac{U_j - U_i}{h}\right)$$

$$\to \mathbb{E}\{[s(X) - \mathbb{E}(s(X)|U)]W(X)\} = \mu.$$

PROOF. $S_{n_2}$ is a $U$-statistic of degree two. Recall $\bar{s}(u) = \mathbb{E}(s(X)|U=u)$. The Hájek projection of $s(X_i)W_j K(\frac{U_j - U_i}{h})$ equals

$$s(X_i)\int_0^1 \bar{W}(v)K\left(\frac{v - U_i}{h}\right)dv + W_j\int_0^1 \bar{s}(u)K\left(\frac{U_j - u}{h}\right)du$$

$$- \int_0^1\int_0^1 \bar{s}(u)\bar{W}(v)K\left(\frac{v-u}{h}\right)dv\,du.$$

Hence, the projection of $S_{n_2}$ equals

$$\hat{S}_{n_2} = n^{-1}\sum_{j=1}^n s(X_j)W_j - \frac{1}{nh}\sum_{i=1}^n s(X_i)\int_0^1 \bar{W}(v)K\left(\frac{v-U_i}{h}\right)dv$$

$$- \frac{1}{nh}\sum_{j=1}^n W_j\int_0^1 \bar{s}(u)K\left(\frac{U_j - u}{h}\right)du$$

$$+ \frac{1}{h}\int_0^1\int_0^1 \bar{s}(u)\bar{W}(v)K\left(\frac{v-u}{h}\right)dv\,du.$$

Furthermore, $\mathbb{E}\{S_{n_2} - \hat{S}_{n_2}\}^2$ is of the order $O(n^{-1}h^{-1}) = o(1)$.

Hence, it remains to show that $\hat{S}_{n_2}$ tends to the desired limit. Now similar to the proof of the previous lemma, it may be shown that

$$\hat{S}_{n_2} - n^{-1} \sum_{j=1}^{n} s(X_j)W_j + n^{-1} \sum_{i=1}^{n} s(X_i)\bar{W}(U_i)$$
$$+ n^{-1} \sum_{j=1}^{n} W_j \bar{s}(U_j) - \int_0^1 \bar{s}(u)\bar{W}(u)\,du \to 0 \quad \text{in probability.}$$

The assertion of the lemma now is a straightforward consequence of the law of large numbers upon noticing that

$$\mathbb{E}[s(X)\bar{W}(U)] = \int_0^1 \bar{s}(u)\bar{W}(u)\,du. \qquad \square$$

The next lemma will be helpful to find the final expansion and limit of $I$.

LEMMA 4.5. *Under the assumptions of Theorem* 2.1, *we have*

$$S_{n_3} \equiv \frac{1}{n(n-1)h}$$
$$\times \sum_{j=1}^{n} \sum_{\substack{i=1 \\ i \neq j}}^{n} s(X_i)W_j \left[ K\left(\frac{F_n(\hat{\beta}^T X_j) - F_n(\hat{\beta}^T X_i)}{h}\right) - K\left(\frac{U_j - U_i}{h}\right) \right]$$
$$\to 0 \quad \text{in probability as } n \to \infty.$$

PROOF. By Taylor's formula,

$$S_{n_3} = \frac{1}{n(n-1)h} \sum_{j=1}^{n} \sum_{\substack{i=1 \\ i \neq j}}^{n} s(X_i)W_j K'(\Delta_{ij})\left[\frac{F_n(\hat{\beta}^T X_j) - F_n(\hat{\beta}^T X_i) - U_j + U_i}{h}\right],$$

where $\Delta_{ij}$ is between the two $K$-ratios in the definition of $S_{n3}$. For each $j$ (and similarly for $i$),

$$|F_n(\hat{\beta}^T X_j) - U_j| \leq |F_n(\hat{\beta}^T X_j) - F^{\hat{\beta}}(\hat{\beta}^T X_j)| + |a_j^{\hat{\beta}}|$$
$$\leq \sup_{t;\theta} |F_n^\theta(t) - F^\theta(t)| + \sup_\theta |a_j^\theta|,$$

where the suprema extend (with large probability) over the set of $\theta$'s with $\|\theta - \beta\| \leq Cn^{-1/2}$. Now it is well known that empirical measures approach the true measure at the rate $O_\mathbb{P}(n^{-1/2})$ uniformly over the class of all halfspaces. See, for example, [26]. In other words, the first supremum is $O_\mathbb{P}(n^{-1/2})$. From Lemma 4.1, the second supremum is $O_\mathbb{P}(n^{-1/2+\alpha})$ uniformly in $1 \leq j \leq n$. Conclude that

$$(4.6) \qquad \sup_{1 \leq j \leq n} |F_n(\hat{\beta}^T X_j) - U_j| = O_\mathbb{P}(n^{-1/2+\alpha}) = O_\mathbb{P}(h).$$

Furthermore, since $K$ has support $[-1, 1]$, the summation in $S_{n_3}$ takes place only w.r.t. those $i, j$ for which at least one of the ratios falls into $[-1, 1]$. If this happens to be true for the first ratio, then by (4.6) also

$$|U_j - U_i| \leq C_3 h,$$

with large probability for some appropriate $C_3$. Summarizing, since $K'$ is bounded, we get, with large probability,

$$|S_{n_3}| \leq C_4 \frac{n^{-1/2+\alpha}}{n(n-1)h^2} \sum_{j=1}^{n} \sum_{\substack{i=1 \\ i \neq j}}^{n} |s(X_i)||W_j| \mathbb{1}_{\{|U_j - U_i| \leq C_3 h\}}.$$

The expectation of the right-hand side is, however, of the order $O(n^{-1/2+\alpha} h^{-1}) = o(1)$. This completes the proof of the lemma. $\square$

We are now ready to analyze the term $I$. From its definition we have

$$n^{-1/2} I = n^{-1} \sum_{j=1}^{n} s(X_j) W_j$$
$$- \frac{1}{n(n-1)h} \sum_{j=1}^{n} \sum_{\substack{i=1 \\ i \neq j}}^{n} s(X_i) W_j K\left(\frac{F_n(\hat{\beta}^T X_j) - F_n(\hat{\beta}^T X_i)}{h}\right).$$

In view of Lemmas 4.4 and 4.5 we therefore get the following result.

COROLLARY 4.6. *Under the assumptions of Theorem 2.1, we have*

$$n^{-1/2} I \to \mu \qquad \text{in probability.}$$

To summarize the results obtained so far, Lemma 4.3 yielded an i.i.d. representation of $n^{-1/2} III$, while Corollary 4.6 provided an in-probability limit for $n^{-1/2} I$. The analysis for $n^{-1/2} IV$ is a bit tricky. At the end it will turn out that it admits an i.i.d. expansion which cancels with the second sum in Lemma 4.3. We may thus conclude that

$$\hat{T}_n = \mu + n^{-1/2} \sum_{i=1}^{n} \varepsilon_i [W_i - \bar{W}(U_i)] + o_\mathbb{P}(1),$$

which coincides with the i.i.d. representation (2.6) of Theorem 2.1. So it remains to show the following representation of $n^{-1/2} IV$.

LEMMA 4.7. *Under the assumptions of Theorem* 2.1,

$$n^{-1/2}IV = n^{-1/2}\sum_{j=1}^{n}\{\Phi(\beta^T X_j)\bar{W}(U_j) - \mathbb{E}[\Phi(\beta^T X_j)\bar{W}(U_j)]\} + o_{\mathbb{P}}(1).$$

PROOF. By Taylor's expansion,

$$-n^{-1/2}IV = \frac{1}{n^{1/2}(n-1)h}$$
$$\times \sum_{\substack{i=1 \\ i\neq j}}^{n}\sum_{j=1}^{n} Y_i^0 W_j \left\{ K\left(\frac{F_n(\hat{\beta}^T X_j) - F_n(\hat{\beta}^T X_i)}{h}\right) - K\left(\frac{U_j - U_i}{h}\right) \right\}$$

(4.7)
$$= \frac{1}{n^{1/2}(n-1)h}$$
$$\times \sum_{\substack{i=1 \\ i\neq j}}^{n}\sum_{j=1}^{n} Y_i^0 W_j K'\left(\frac{U_j - U_i}{h}\right) \frac{F_n(\hat{\beta}^T X_j) - F_n(\hat{\beta}^T X_i) - U_j + U_i}{h}$$
$$+ \frac{1}{2n^{1/2}(n-1)h} \sum_{\substack{i=1 \\ i\neq j}}^{n}\sum_{j=1}^{n} Y_i^0 W_j K''(\Delta_{ij}) \frac{[\cdots]^2}{h^2},$$

where $\Delta_{ij}$ is between the two $K$-ratios in the representation of $IV$. We shall show that the second double sum is negligible, while the first contributes to the i.i.d. representation of $\hat{T}_n$. First, we write

$$F_n(\hat{\beta}^T X_j) = F_n(\hat{\beta}^T X_j) - F^{\hat{\beta}}(\hat{\beta}^T X_j) - F_n^{\beta}(\beta^T X_j) + F^{\beta}(\beta^T X_j)$$
$$+ F^{\hat{\beta}}(\hat{\beta}^T X_j) + F_n^{\beta}(\beta^T X_j) - F^{\beta}(\beta^T X_j),$$

and similarly for the index $i$. The first line equals, with $\theta = \hat{\beta}$ and $x = \hat{\beta}^T X_j, y = \beta^T X_j$, the quantity $G_n^{\theta}(x,y)$ appearing in Lemma 4.2. Conclude from that result that

(4.8)
$$\frac{1}{n^{1/2}(n-1)h^2}\sum_{i=1}^{n}\sum_{j=1}^{n}|Y_i^0 W_j|\left|K'\left(\frac{U_j - U_i}{h}\right)\right||G_n^{\theta}(x,y)|$$
$$= O_{\mathbb{P}}(\sqrt{n^{-1/2+\alpha}\ln n})\frac{1}{n^2 h^2}\sum_{i=1}^{n}\sum_{j=1}^{n}|Y_i^0 W_j|\left|K'\left(\frac{U_j - U_i}{h}\right)\right|.$$

The double sum is easily seen to be bounded in probability. Since

$$n^{-1/2+\alpha}\ln n \to 0,$$

this proves that (4.8) tends to zero in probability. Next we study

$$\frac{1}{n^{1/2}(n-1)h^2}\sum_{i=1}^{n}\sum_{j=1}^{n}Y_i^0 W_j K'\left(\frac{U_j-U_i}{h}\right)$$
$$\times [F_n^\beta(\beta^T X_j)-F^\beta(\beta^T X_j)-F_n^\beta(\beta^T X_i)+F^\beta(\beta^T X_i)].$$

This sum is a $V$-statistic (see [27]), with a kernel depending on $h$ and hence on $n$. It is asymptotically equal to a $U$-statistic whose Hájek projection equals

$$h^{-2}\int_0^1\int_0^1 \psi(u)\bar{W}(v)K'\left(\frac{v-u}{h}\right)[\bar{\alpha}_n(v)-\bar{\alpha}_n(u)]\,du\,dv.$$

Here, $\bar{\alpha}_n$ is the (uniform) empirical process pertaining to the $U_j$'s. Transformation of integrals, $C$-tightness of $\bar{\alpha}_n$, $n\geq 1$, and the fact that $K'$ has compact support $[-1,1]$ yield that the last double integral is equivalent to

$$h^{-1}\int_0^1\int_{-1}^1 \psi(v-wh)\bar{W}(v)K'(w)[\bar{\alpha}_n(v)-\bar{\alpha}_n(v-wh)]\,dw\,dv.$$

By continuity of $\psi$, this is asymptotically equivalent to

(4.9)
$$h^{-1}\int_0^1\int_{-1}^1 \psi(v)\bar{W}(v)K'(w)[\bar{\alpha}_n(v)-\bar{\alpha}_n(v-wh)]\,dw\,dv$$
$$=-h^{-1}\int_0^1\int_{-1}^1 \psi(v)\bar{W}(v)K'(w)\bar{\alpha}_n(v-wh)\,dw\,dv.$$

Check that

$$\frac{1}{h}\int_{-1}^1 K'(w)\bar{\alpha}_n(v-wh)\,dw = \sqrt{n}[\bar{f}_n(v)-1].$$

Here

$$\bar{f}_n(v) = \frac{1}{nh}\sum_{i=1}^n K\left(\frac{v-U_i}{h}\right)$$

is the kernel density estimator for the uniform sample $U_1,\ldots,U_n$.

Hence, (4.9) equals

(4.10) $$-\sqrt{n}\int_0^1 \psi(v)\bar{W}(v)[\bar{f}_n(v)-1]\,dv.$$

Introducing the smoothed empirical distribution,

$$d\tilde{F}_n = \bar{f}_n\,dv,$$

and the pertaining empirical process $\tilde{\alpha}_n = \sqrt{n}(\tilde{F}_n - Id)$, where $Id$ denotes the identity function on $(0,1)$, (4.10) becomes

$$-\int_0^1 \psi(v)\bar{W}(v)\tilde{\alpha}_n(dv).$$

It is known that

$$(4.11) \qquad \int_0^1 \psi \bar{W} \, d\tilde{\alpha}_n = \int_0^1 \psi \bar{W} \, d\bar{\alpha}_n + o_{\mathbb{P}}(1).$$

A simple proof of (4.11) may be obtained by using oscillation results for empirical processes; see [30]. We shall shortly see that all other terms will be negligible for the i.i.d. representation of $n^{-1/2} IV$, so that

$$(4.12) \qquad n^{-1/2} IV = \int_0^1 \psi \bar{W} \, d\bar{\alpha}_n + o_{\mathbb{P}}(1),$$

as desired. To justify (4.12), we next bound

$$(4.13) \quad \frac{1}{n^{1/2}(n-1)h^2} \sum_{i=1}^n \sum_{j=1}^n Y_i^0 W_j K'\!\left(\frac{U_j - U_i}{h}\right)(\hat{U}_j - U_j - \hat{U}_i + U_i),$$

where

$$\hat{U}_j = F^{\hat{\beta}}(\hat{\beta}^T X_j), \qquad 1 \le j \le n.$$

Hence, the $\hat{U}_j$ and $U_j$ incorporate the theoretical distribution functions

$$F^\theta(x) = \mathbb{P}(\theta^T X \le x)$$

at $\theta = \hat{\beta}$ and $\theta = \beta$, respectively. From Lemma 4.1,

$$(4.14) \qquad \max_{1 \le j \le n} |\hat{U}_j - U_j| = O_{\mathbb{P}}(n^{-1/2+\alpha}).$$

This bound will sometimes be helpful to further simplify (4.13). First, because $K'$ is an odd function, (4.13) may be written as

$$(4.15) \quad \frac{1}{n^{1/2}(n-1)h^2} \sum_{i=1}^n \sum_{j=1}^n (Y_i^0 W_j + Y_j^0 W_i) K'\!\left(\frac{U_j - U_i}{h}\right)(\hat{U}_j - U_j).$$

We shall only deal with the sum involving $Y_i^0 W_j$, the other being dealt with in a similar way. Now

$$\frac{1}{n^{1/2}(n-1)h^2} \sum_{j=1}^n W_j (\hat{U}_j - U_j) \sum_{i=1}^n Y_i^0 K'\!\left(\frac{U_j - U_i}{h}\right)$$

$$= \frac{1}{(n-1)h} \sum_{j=1}^n W_j (\hat{U}_j - U_j) \left[\frac{1}{\sqrt{n}h} \sum_{i=1}^n \varepsilon_i K'\!\left(\frac{U_j - U_i}{h}\right)\right]$$

$$+ \frac{1}{(n-1)h} \sum_{j=1}^n W_j (\hat{U}_j - U_j)$$

$$\times \left[\frac{1}{\sqrt{n}h} \sum_{i=1}^n \left\{\Phi(\beta^T X_i) K'\!\left(\frac{U_j - U_i}{h}\right) - \mathbb{E}[\cdots]\right\}\right]$$

$$+ \frac{1}{\sqrt{n}h^2} \sum_{j=1}^n W_j(\hat{U}_j - U_j) \int_0^1 \psi(v) K'\left(\frac{U_j - v}{h}\right) dv.$$

In the first two double series, first apply (4.14) to bound $|\hat{U}_j - U_j|$ uniformly in $j$. The expectation of, for example,

$$\frac{1}{n} \sum_{j=1}^n |W_j| \left| \frac{1}{\sqrt{n}h} \sum_{i=1}^n \varepsilon_i K'\left(\frac{U_j - U_i}{h}\right) \right|,$$

is easily seen to be bounded. Similarly for the second series. Conclude that each sum is

$$O_{\mathbb{P}}(h^{-1} n^{-1/2 + \alpha}) = o_{\mathbb{P}}(1).$$

As to the last $j$-sum, substitute $w = \frac{U_j - v}{h}$, apply Taylor's expansion to $\psi(U_j - wh)$ and use the fact that

$$\int_{-1}^1 K'(w)\, dw = 0, \qquad \int_{-1}^1 w K'(w)\, dw = -1$$

to finally get that the last sum equals

$$\frac{1}{\sqrt{n}} \sum_{j=1}^n W_j \psi'(U_j)(\hat{U}_j - U_j) + o_{\mathbb{P}}(1)$$

$$= \frac{1}{\sqrt{n}} \sum_{j=1}^n \bar{W}(U_j) \psi'(U_j)(\hat{U}_j - U_j) + o_{\mathbb{P}}(1).$$

Similar arguments yield for the double sum in (4.15) including the factors $Y_j^0 W_i$ the representation

$$\frac{1}{\sqrt{n}} \sum_{j=1}^n \bar{W}'(U_j) \psi(U_j)(\hat{U}_j - U_j) + o_{\mathbb{P}}(1).$$

Conclude that so far we have shown that (4.13) equals

(4.16) $$n^{-1/2} \sum_{j=1}^n (\bar{W}\psi)'(U_j)(\hat{U}_j - U_j) + o_{\mathbb{P}}(1).$$

At this point we see that another simple application of (4.14) even for bounded $X$'s, that is, $\alpha = 0$, does not yield an $o_{\mathbb{P}}(1)$ term. Therefore, we have to analyze $\hat{U}_j$ and $U_j$ in a different way. As we shall see, finally, and in a disguised form, we take advantage of the fact that, for each $\theta$ every projection $\theta^T X$ of $X$ is transformed into a uniform random variable $F^\theta(\theta^T X)$. Fix such a $\theta$ and note that, for a random vector $X$ with the same distribution as $X_1$ but being independent of the sample $(X_i, Y_i), 1 \leq i \leq n$, one gets

$$F^\theta(\theta^T X_j) = \mathbb{E}\{\mathbb{1}_{\{\theta^T X \leq \theta^T X_j\}} | \mathcal{F}_n\}.$$

Here $\mathcal{F}_n = \sigma(X_i, Y_i, 1 \leq i \leq n)$ is the $\sigma$-field generated by the observations. Conclude that, for $\theta = \hat{\beta}$,

$$\hat{U}_j - U_j = \mathbb{E}\{\mathbb{1}_{\{\theta^T X \leq \theta^T X_j\}} - \mathbb{1}_{\{\beta^T X \leq \beta^T X_j\}}|\mathcal{F}_n\}$$
$$= \mathbb{E}\{\mathbb{1}_{\{\theta^T X \leq \theta^T X_j\}} - \mathbb{1}_{\{\theta^T X \leq \beta^T X_j\}}|\mathcal{F}_n\}$$
$$+ \mathbb{E}\{\mathbb{1}_{\{\theta^T X \leq \beta^T X_j\}} - \mathbb{1}_{\{\beta^T X \leq \beta^T X_j\}}|\mathcal{F}_n\},$$

whence

$$n^{-1/2} \sum_{j=1}^{n} (\bar{W}\psi)'(U_j)(\hat{U}_j - U_j)$$

$$(4.17) \quad = n^{-1/2} \sum_{j=1}^{n} (\bar{W}\psi)'(U_j)(\hat{\beta} - \beta)^T X_j f^\beta(\beta^T X_j) + o_\mathbb{P}(1)$$

$$(4.18) \quad + \mathbb{E}\left\{ n^{-1/2} \sum_{j=1}^{n} (\bar{W}\psi)'(U_j)[\mathbb{1}_{\{\theta^T X \leq \beta^T X_j\}} - \mathbb{1}_{\{\beta^T X \leq \beta^T X_j\}}]\Big|\mathcal{F}_n \right\}.$$

The process inside the conditional expectation is, after centering, asymptotically $C$-tight. With $\theta = \hat{\beta} \to \beta$, we therefore obtain

$$\mathbb{E}\{\cdots|\mathcal{F}_n\} = n^{1/2} \mathbb{E}\left\{ \int_{F^\beta(\theta^T X)}^{F^\beta(\beta^T X)} (\bar{W}\psi)'(u)\, du \Big| \mathcal{F}_n \right\} + o_\mathbb{P}(1)$$
$$= n^{1/2}(\beta - \hat{\beta})\mathbb{E}\{(\bar{W}\psi)'(U)X f^\beta(\beta^T X)\} + o_\mathbb{P}(1),$$

where the last equality follows from the mean value theorem, $n^{1/2}(\hat{\beta} - \beta) = O_\mathbb{P}(1)$ and the facts that $\hat{\beta}$ is measurable with respect to $\mathcal{F}_n$ and $X$ is independent of $\mathcal{F}_n$. Inserting this in (4.17) and (4.18), we thus get

$$n^{-1/2} \sum_{j=1}^{n} (\bar{W}\psi)'(U_j)(\hat{U}_j - U_j)$$
$$= n^{1/2}(\hat{\beta} - \beta) n^{-1} \sum_{j=1}^{n} \{(\bar{W}\psi)'(U_j) X_j f^\beta(\beta^T X_j) - \mathbb{E}[\cdots]\} + o_\mathbb{P}(1).$$

Since $n^{1/2}(\hat{\beta} - \beta)$ is stochastically bounded and the sample mean tends to zero according to the SLLN, this shows that (4.16) tends to zero in probability.

It remains to bound (4.7), but this is easy. In view of Lemma 4.2, upon applying by now standard arguments, we have

$$|(4.7)| = o_\mathbb{P}(1).$$

This completes the proof of Lemma 4.7. □

**Acknowledgments.** We are grateful to the referees and an Associate Editor who pointed out to us relevant literature missing in an earlier version of the paper.

Mathematical Institute
Universität Giessen
Arndtstrasse 2
D-35392 Giessen
Germany
e-mail: Winfried.Stute@math.uni-giessen.de

Department of Statistics
and Actuarial Science
The University of Hong Kong
Pokfulam Road, Hong Kong
China
e-mail: lzhu@hku.hk